\documentclass[12pt]{amsart}
\usepackage{amsmath,amssymb}

\newtheorem{thm}{Theorem}[section]
\newtheorem{prop}[thm]{Proposition}

\newtheorem{cor}[thm]{Corollary}
\newtheorem{ex}[thm]{Example}

\newcommand\sub{\subseteq}
\newcommand\imp{\Rightarrow}

\begin{document}


    \title[Maximal divisorial ideals and $t$-maximal ideals]
{Maximal divisorial ideals and \\ $t$-maximal ideals}



\author{Stefania Gabelli}

\address{Dipartimento di Matematica, Universit\`{a} degli Studi Roma Tre,
Largo S.  L.  Murialdo,
1, 00146 Roma, Italy}

\email{gabelli@mat.uniroma3.it}

\author{Moshe Roitman}

\address{Department of Mathematics, University of
Haifa, Mount Carmel, Haifa 31905, Israel}

\email{mroitman@math.haifa.ac.il}




\begin{abstract}
We give conditions for a maximal divisorial ideal to be $t$-maximal and
show with examples that, even in a completely integrally closed domain,
maximal divisorial ideals need not be $t$-maximal.
\end{abstract}

    \maketitle




\section*{Introduction}

The $v$-operation and the $t$-operation are the the best known and most 
useful  star operations;
mainly because the 
structure of certain semigroups of $t$-ideals reflects the 
multiplicative properties of an integral domain. In this context an 
important role is played by the 
prime and the maximal $v$- and $t$-ideals.

Since the $t$-operation is a star operation of finite type, a domain $R$
has always $t$-maximal ideals. On the other hand, the set of $v$-maximal
ideals may be empty.

In this paper we deal with the following question:

\smallskip
{\it Assume that $M$ is a $v$-maximal ideal of $R$, is $M$ necessarily 
a $t$-maximal ideal?}
\smallskip

We show that although the answer is positive in a large class of domains, namely 
in the class of $v$-coherent domains, it is negative in general. 
In fact we give two examples of a
$v$-maximal ideal $P$ that is not a
$t$-maximal ideal.  In the first example $P$ is an upper to zero of 
a completely integrally closed polynomial ring, thus
$P$ is $v$-invertible. In the second 
example $P$ is a strongly divisorial ideal of an integrally closed 
semigroup ring.




\section{Preliminaries and notations}

Throughout this paper $R$ will denote an integral domain
with quotient field $K$. We will refer to a fractional ideal 
as an {\it ideal} and will call a fractional ideal contained in $R$ an 
{\it integral ideal}.

We recall that a star operation is an application $I \to I^*$ 
from the set $F(R)$ of nonzero ideals of $R$ to itself such that: 
\begin{enumerate}
	\item $R^* = R$ and $(aI)^* = aI^*$, for all $a \in K \smallsetminus \{0\}$;
	\item $I \sub I^*$ and $I \sub J \imp I^* \sub J^*$;
	\item $I^{**} = I^*$.
\end{enumerate}

General references for systems
of ideals and star operations are 
 \cite{ G, Gr, HK, J}.

We denote by $f(R)$ the set of nonzero finitely generated ideals of $R$.
A star operation $*$ is of {\it finite type} if, for each $I \in F(R)$, 
 $I^* = \cup \{J^*\, |\, J \sub I$ and $J \in f(R) \}$.
To any  star operation $*$, we can associate a star operation
$*_{f}$ of finite type
by defining $I^{*_{f}}= \cup \{J^*\, |\, J \sub I$ and $J \in f(R) \}.$ 
Clearly $I^{*_{f}}\sub I^*$.

The $v$- and the
$t$-operations are particular star operations, defined in the following 
way.

For a pair of nonzero ideals $I$
and $J$ of a domain $R$ we let $(J:I)$ denote the set $\{x\in K|\,
xI\subseteq J\}$ and $(J:_{R}I)$ denote the set $\{x\in R|\, xI\subseteq
J\}$. 
We set $I_v=(R:(R:I))$ and $I_t=\bigcup J_v$ with the
union taken over all finitely generated ideals $J$ contained
in $I$.

The $t$-operation is the finite type 
star operation associated to the $v$-operation.

 A nonzero ideal $I$ is called a $*$-{\it ideal} if $I= I^*$.
 Thus a nonzero  ideal $I$ is a {\it $v$-ideal}, or is {\it divisorial}, if
$I=I_v$, and it is a {\it $t$-ideal} if $I=I_t$. Note that $I$ is a 
$t$-ideal if and only if $J_{v}\sub I$ whenever $J$ is finitely generated 
and $J\sub I$.
 
 The  set $F_{*}(R)$ of $*$-ideals of 
 $R$ is a semigroup with respect to the 
 $*$-{\it multiplication}, defined by $(I,J) \to (IJ)^*$, with unity $R$. 
 
 We say that an ideal $I \in F(R)$ is $*$-{\it invertible} if $I^*$ is a 
 unit in the 
 semigroup $F_{*}(R)$. In this case the $*$-{\it inverse} of 
$I$ is $(R:I)$.
 Thus $I$ is  $*$-invertible if and only if $(I(R:I))^*=R$. Invertible 
 ideals are ($*$-invertible) $*$-ideals.

 We have $I \sub I^* \sub I_{v}$, so that any divisorial ideal is a 
$*$-ideal and any $*$-invertible  ideal is $v$-invertible.
In particular a divisorial  ideal is a $t$-ideal and a $t$-invertible 
ideal is $v$-invertible. 

 A nonzero ideal $I$ is $*$-{\it finite} if $ I^* =J^*$ for 
 some finitely generated ideal $J$.
Since the $v$- and the $t$-operation coincide on finitely generated ideals
and since $I_t=J_t$ implies $I_v=J_v$,
 an ideal $I$ is $t$-finite if and only if $I_{v}= J_{v}$
 (equivalently $(R:I) = (R:J)$)  for 
 some finitely generated ideal $J \sub I$.
It follows that the set $f_{v}(R)$ of the $v$-finite divisorial  
ideals coincides with the set of the $t$-finite  $t$-ideals. $f_{v}(R)$ 
 is a sub-semigroup of $F_{v}(R)$.

An ideal $I$ is $t$-invertible if 
and only if it is $v$-invertible and both $I$ and $(R:I)$ are $t$-finite. 
Hence the set of the $t$-invertible $t$-ideals of $R$, here denoted by  
$T(R)$,
is the largest subgroup of $f_{v}(R)$. The importance of the notion of 
$t$-invertibility is well illustrated in \cite{Z}.
 
Denoting by $Inv(R)$ the group of the invertible ideals of $R$, we have 
$$Inv(R) \sub T(R) \sub f_{v}(R) \sub F_{v}(R) \sub F_{t}(R) \sub F(R)$$ 
and 
$$Inv(R) \sub f(R).$$ 
Several important classes of domains may be characterized by the fact that some 
of these inclusions are
equalities.  For example $R$ is a {\it 
Pr\"ufer} domain if and only if  $Inv(R) = f(R)$ \cite{G}, 
it is a {\it Krull} domain if and only if  
$T(R)=F_{t}(R)$ \cite{K2} and it is a {\it 
Pr\"ufer $v$-multiplication} domain, for short a P$v$MD, if and only if  
$T(R) =  f_{v}(R)$ \cite{Gr} . A {\it Mori} domain 
is a domain satisfying the ascending chain condition on integral
divisorial ideals and has the property that  $f_{v}(R)=F_{t}(R)$.
Noetherian and Krull domains are Mori.  A recent 
reference for Mori domains is \cite{Bar}. 
The class of domains 
with the property  that $F_{v}(R)=F(R)$ have been studied by several authors \cite{B, He, 
M, Ba, BS}.
A domain such that $F_{v}(R) = F_{t}(R)$ is called in \cite{HZ} a $TV$-domain. Examples 
of $TV$-domains are given in \cite{HZ, K1, HC}. Mori and  
pseudovaluations domains which are not valuation domains
are $TV$-domains.




\section{When a maximal divisorial ideal is $t$-maximal}

A prime $*$-ideal is  called a $*$-prime. A $*$-{\it maximal} ideal is an ideal
that is maximal in the set of the proper integral
$*$-ideals. A $v$-maximal ideal is also called a {\it maximal divisorial 
ideal}.
A $*$-maximal ideal is a prime ideal (if it exists). 

If $*$ is a star operation of finite type, an easy application of the
Zorn Lemma shows that the set $*$-Max$(R)$ of the $*$-maximal ideals of
$R$ is not empty. Moreover, for each $I\in F(R)$, $I^* = \bigcap_{M\in
*\text{-Max}(R)}I^*R_{M} $ \cite{Gr}. In particular the set of the
$t$-maximal ideals is not empty and $I_{t} = \bigcap_{M\in 
t\text{-Max}(R)}I_{t}R_{M} $. On the contrary, the set of maximal divisorial
ideals may be empty, like for example when $R$ is a rank-one nondiscrete
valuation domain.

If $M$ is a $*$-maximal ideal that is not $*$-invertible, then 
$M=(M(R:M))^{*}$ and so
$(M:M)=(R:M)$. An ideal $I$ with the property that $(R:I)=(I:I)$ is called {\it 
strong}. 
A strong ideal is never $*$-invertible and we have just seen 
that a $*$-maximal ideal is either $*$-invertible or strong. 

An ideal which is strong and divisorial is called {\it strongly 
divisorial}.

\begin{prop}\label{max}
	
	\

	\begin{enumerate}
	
\item If $M$ is a maximal divisorial ideal of $R$, then $M = x^{-1}R\cap 
R$, for some element $x \in K$. Hence $(R:M) = (R+xR)_{v}$.

\item   If $P$ is a prime divisorial ideal of $R$ such
that $(R:P)=R+xR$, for some element $x \in K$,
then $P$ is maximal divisorial.  
\end {enumerate}
\end{prop}

\begin{proof} 
	(1) If $x\in (R:M)\smallsetminus R$, then $M \sub x^{-1}R$ and $R 
	\nsubseteq x^{-1}R$. 
	Since an intersection of divisorial ideals 
is divisorial and $M$ is $v$-maximal, we have $M = x^{-1}R\cap 
R = (R:R+xR)$.
	
	(2) Let $Q$ be a proper divisorial ideal containing $P$.
Since $Q$ is divisorial, $(R:Q)\nsubseteq R$.
Since $(R:Q)\subseteq (R:P)=R+xR$, we see that
there exists an element $y\in R$ such that $xy\in (R:Q)\smallsetminus R$. Thus
 $y\notin P$, and $xyQ\subseteq R$.  Since $P=(R:R+xR)$, we obtain that
$yQ\subseteq P$.  Since $P$ is a
prime ideal, we conclude that $Q\subseteq P$. Hence $P$ is
maximal divisorial.  
\end{proof}

In a Mori domain $R$, all the prime divisorial ideals are of the form $x^{-1}R\cap 
R = (R:R+xR)$ \cite[Corollary 2.5]{HLV}.

A domain has the property that each $t$-maximal ideal is divisorial 
if and only if every ideal $I$ such that $(R:I) = R$ is $t$-finite 
\cite[Proposition 2.4]{HZ}. A domain of this type is called an $H$-domain 
in \cite{GV}. A $TV$-domain is clearly an $H$-domain, but the converse 
is not true \cite{HZ, BGR}.

The following proposition gives conditions for a  divisorial prime ideal 
to be a $t$-maximal ideal. A proof can be found in \cite{Ga}.

\begin{prop}\label{vtmax}

\
\begin{enumerate}
	\item A $v$-invertible  divisorial prime is maximal divisorial;
	\item A $v$-finite maximal divisorial ideal is $t$-maximal;
	\item A $v$-finite $v$-invertible  divisorial prime is $t$-invertible;
	\item A $t$-invertible $t$-prime is $t$-maximal.
	\end{enumerate}
	\end{prop}

	We remark that in general a $*$-invertible $*$-prime need not be 
	$*$-maximal (for example a principal prime ideal is not necessarily a 
	maximal ideal) and that a (non-prime) $v$-finite $v$-invertible 
	divisorial ideal need not be $t$-invertible \cite{D}.
	
\begin{cor}\label{cormax}
Assume that each maximal divisorial ideal of $R$ is a $t$-maximal ideal. 
Then each $v$-invertible divisorial prime is a $t$-invertible $t$-maximal ideal.
\end{cor}
\begin{proof}
Let $P$ be a $v$-invertible divisorial prime. By Proposition \ref{vtmax}, 
$P$ is maximal divisorial and so $t$-maximal. Since $P$ is not strong, 
then it is  $t$-invertible.
\end{proof}

In general, if each $v$-invertible divisorial prime of $R$ is $t$-invertible, 
it is not true that each  $v$-invertible ideal is $t$-invertible. This 
last property is in fact equivalent to $R$ being an $H$-domain \cite
 [Proposition 4.2]{Z}. The ring of entire functions is not an $H$-domain, but
all its divisorial primes are $t$-invertible (see for example \cite[Section 
2]{Ga}).
	
	A $v$-{\it coherent} domain is a domain $R$ with the property that,
	for each finitely generated ideal $J$, the ideal $(R:J)$ is $v$-finite.
	This class of domains was first studied (under a different name) in 
	\cite{elAb} and is very large, properly including 
	P$v$MD's, Mori domains and coherent domains \cite{elAb, GH}.
	(A domain is {\it 
	coherent} if each finitely generated 
	ideal is finitely presented, or, equivalently, if  the intersection of each pair of finitely generated ideals 
	is finitely generated.)

\begin{prop}
	If $R$ is $v$-coherent, then each maximal divisorial ideal is $t$-maximal.
\end{prop}
\begin{proof}
	Let $M$ be a maximal divisorial ideal of $R$. Then $M=x^{-1}R\cap R=(R:R+xR)$ for some 
	$x\in K$ (Proposition \ref{max}). Since $R$ is $v$-coherent, then $M$ is
	$v$-finite and so $t$-maximal by Proposition \ref{vtmax}.
	\end{proof}

A domain $R$ is {\it completely integrally closed} if and only if 
        $F_{v}(R)$ is a group under $v$-multiplication \cite{G}. 
        If $F_{v}(R) = T(R)$, then $R$ is a completely
integrally closed $H$-domain, equivalently a Krull domain \cite{Ga, GV}.

A  divisorial prime of a completely integrally closed domain, being 
$v$-invertible, is always maximal divisorial by Proposition \ref{vtmax}. 
We will see in the next section that it 
need not be $t$-maximal. As a matter of fact, a divisorial 
prime $P$ of a completely integrally 
closed domain has height one and $P$ is $t$-maximal if and only if it is 
$v$-finite, if and only if it is $t$-invertible \cite[Theorem 2.3]{Ga}. 

A completely integrally closed $v$-coherent domain is a (completely integrally closed)
P$v$MD. In this case each divisorial prime is $t$-maximal by 
Corollary \ref{cormax}.

\bigskip


We now turn to the case of polynomial rings.

	We denote by $\mathbf X$ 
a set of independent indeterminates over $R$ and by $R[\mathbf X]$ the 
polynomial ring in this set of indeterminates.

It is well known that the correspondence $I \mapsto I[\mathbf X]$ induces  
inclusion preserving injective maps
$t(R) \longrightarrow t(R[\mathbf X])$ and 
$D(R) \longrightarrow D(R[\mathbf X])$. Moreover, 
$M$ is a $t$-maximal ideal, respectively a maximal 
divisorial ideal, of $R[\mathbf X]$ such that $M\cap 
	R \neq (0)$, if and only if $M=(M\cap R)[\mathbf X]$ and  $M\cap 
	R$ is a $t$-maximal ideal, respectively a maximal 
divisorial ideal, of $R$ (see for example \cite[Lemma 2.1]{FGH} and 
\cite[Theorem 3.6]{Rt}). 

Thus, if
each maximal divisorial ideal 
	of $R[\mathbf X]$ is $t$-maximal, $R$ has the same property.

	On the other hand, Example 3.1 in the next section shows that if $M \cap R = 
	(0)$, then $M$ may be maximal divisorial but not $t$-maximal.
	
		A prime ideal $Q$ of $R[\mathbf X]$ such that $Q \cap R = 
	(0)$ is called an {\it upper to zero}. $Q$ is an upper to zero of
	height one if and only if $Q = fK[\mathbf X]\cap R[\mathbf X]$ 
	for some polynomial $f \in R[\mathbf X]$, irreducible in $K[\mathbf X]$
	\cite[Lemma 2.1]{GHL}. In one indeterminate, all 
	the uppers to zero are of this form.
	
Recall that if $R$ is integrally closed and $f$ is a nonzero polynomial
of $R[\mathbf X]$, 
then $fK[\mathbf X]\cap R[\mathbf X] = f(R:c(f))[\mathbf X]$ 
\cite[Corollary 34.9]{G}.
(Here $c(f)$ denotes the {\it content} of 
$f$, that is the fractional ideal of $R$ generated by the coefficients 
of $f$.) Hence if $R$ is integrally closed, an upper to zero of height one   
is always divisorial and if $R$ is completely integrally closed, an upper to 
zero of height one, being $v$-invertible, is always maximal divisorial.

In general, an upper to zero 
is $t$-maximal if and only if it is $t$-invertible; in this case it
has height one \cite[Section 3]{GHL}. We now show that a similar result holds for the 
$v$-operation.

\begin{prop}\label{upperdiv}
A divisorial upper to zero is a maximal divisorial ideal if and only if 
it is $v$-invertible. In this case it has height one.
	\end{prop}
	\begin{proof}
		A divisorial $v$-invertible prime is always maximal divisorial 
		(Proposition \ref{vtmax} (1)).
		
		Conversely, let $P \subseteq R[\mathbf X]$ be an upper to zero that is
		maximal divisorial. Then
		$P = \frac fg R[\mathbf X]\cap R[\mathbf X] \subseteq 
		fK[\mathbf X]\cap R[\mathbf X]$ , for some $f, g \in R[\mathbf X]$, 
		$g \neq 0$ (Proposition \ref{max}(1)). Since $P\cap R = (0)$ 
		and $f=\frac fg g \in P$, then $f\notin R$. We may also assume that 
		$f$ and $g$ are coprime in $K[\mathbf X]$.
		
		Let $h = f\alpha \in fK[\mathbf X]\cap R[\mathbf 
		X]$, with $\alpha \in K[\mathbf X]$. 
		There is a 
		nonzero $c\in R$ such that $c\alpha \in R[\mathbf X]$. Hence
		$ch = (c\alpha)f = (c\alpha g)\frac fg \in 
		P$. Since $c \notin P$, then $h \in P$. 
		 
		We conclude that 
		$P = fK[\mathbf X]\cap R[\mathbf X]$ has height one.
		
		In addition, $\frac gf \in (R[\mathbf X] : P)$, but $\frac gf 
		\notin (P:P)$. Otherwise $g = \frac gf f \in P$ and so $g = \frac 
		fg t$ for some $t\in R[\mathbf X]$. Then $f$ divides $g^2$ in $K[\mathbf 
		X]$, which is impossible, because $f$ and $g$ are coprime and $f 
		\notin K$. 
	
		It follows that $P$ is not strong and, being maximal divisorial, 
		is $v$-invertible.
		\end{proof}

		The following result was proved in \cite{GV} for one indeterminate.
		
	\begin{prop}\label{Hpoly} $R$ is an $H$-domain if and only if $R[\mathbf X]$
		is an $H$-domain.
		\end{prop}
		\begin{proof} An extended prime $P[\mathbf X]$ is a $t$-maximal ideal, 
		respectively a maximal divisorial ideal, if and only if so is $P$, 
		\cite[Lemma 2.1]{FGH} and \cite[Theorem 3.6]{Rt}.
	A $t$-maximal upper to zero is $t$-invertible 
	by \cite[Theorem 2.3]{GHL}. Hence it is divisorial. 
	\end{proof}
		 
		The domain $R$ is 
said to be a $UMT$-{\it domain} if every upper to zero of $R[X]$ 
is a $t$-maximal 
ideal \cite{HZ2}. This property is stable under polynomial extensions, 
in fact $R$ is a $UMT$-domain if and only if $R[\mathbf X]$ is a 
$UMT$-domain \cite [Theorem 2.4]{FGH}.
The integrally closed $UMT$-domains are exactly the 
P$v$MDs \cite[Proposition 3.2]{HZ2}.

The following proposition is immediate.

\begin{prop}\label{poly}
	Assume that $R$ is an $UMT$-domain. Then each maximal divisorial ideal 
	of $R$ 	is $t$-maximal if and only if $R[\mathbf X]$ has the same property.
	\end{prop}
	
	We conclude this section recalling that it is not known whether $R$ 
	$v$-coherent implies that $R[X]$ is $v$-coherent. This is true under the 
	additional hypothesis that $R$ is integrally closed \cite{elAb}. In 
	this case, each prime of $R[X]$ upper to zero is divisorial $v$-finite. 
	When $R$ is $v$-coherent and completely integrally closed (thus a 
	completely integrally closed P$v$MD), each upper to 
	zero of $R[X]$ is $t$-maximal (and $t$-invertible).



\section{maximal divisorial ideals that are not $t$-maximal}

In this section we give two examples of a
maximal divisorial ideal $P$ of an integral domain $R$ that is not a
$t$-maximal ideal.  In the first example $R$ is a completely integrally closed
polynomial ring in one indeterminate and $P$ is an upper to zero, thus
$P$ is $v$-invertible. In the second 
example $R$ is an integrally closed semigroup ring and $P$ is strongly divisorial.


\begin{ex}
An upper to zero $P$ of a completely integrally closed polynomial ring $R[X]$ 
that is maximal divisorial but not t-maximal. $P$ is necessarily v-invertible.
\end{ex}

{\it
Let $y,z$ and $\mathbf t=\{t_n (n\ge1)\}$ be independent indeterminates over a
field $k$. Let $S$ be the semigroup of monomials $f$ of $k[y, z, \mathbf t]$ 
satisfying the conditions 
$\deg_{y,z} f\ge \deg_{ t_n}f$ for all $n\ge 1$, and let $R = k[S]$ the 
semigroup ring over $k$ generated by $S$. 

Set 
$$P=(y+zX)K[X]\cap R[X],$$
where $K$ is the field of fractions of $R$ 
and $X$ is an indeterminate over $R$.
Then $R$ (and so also $R[X]$) is completely integrally closed, and $P$ is a maximal
divisorial ideal of $R[X]$ that is not t-maximal.}

\begin{proof}
	
	\

\begin{enumerate}
\item
{\em $R[X]$ is completely integrally closed.}

It is enough to show that $R$ is completely integrally closed. Since
$R=k[S]$ is a semigroup ring over the field $k$, by \cite[Corollary 12.7
(2)]{gilmer} to this end it suffices to show that the semigroup $S$ is
completely integrally closed.

Let $u,v,w\in S$ so that $u(\frac v w)^m\in S$ for all $m\ge 1$. Fix
$n\ge 1$. Then $\deg_{y,z}(u(\frac vw)^m)\ge \deg_{t_n}(u(\frac vw)^m)$
for all $m$. Hence $$\deg_{y,z} u+m\deg_{y,z} (\frac vw)\ge\deg_{
t_n}u+m\deg_{t_n} (\frac vw).$$ Divide by $m$ and let $m$ go to $\infty$
to obtain that $\deg_{y,z}(\frac vw)\ge \deg_{t_n}(\frac vw)$. The same
argument shows that $\frac vw$ is a monomial, that is has a nonnegative
degree in each indeterminate. It follows that $\frac vw\in S$; thus $S$
is completely integrally closed.

\item
{\em $P$ is an upper to zero of $R[X]$ that is a $v$-invertible maximal divisorial 
ideal.}

$P$ is clearly an upper to zero. Since $R$ is integrally closed, then
$P=(y+zX)(R:(y,z))[X]$ by \cite[Corollary 34.9]{G}, hence $P$
is divisorial. But $R[X]$ is completely integrally closed; thus 
$P$ is $v$-invertible and so maximal divisorial (Proposition \ref{vtmax}).

\item
{\em $P$ is not t-maximal.}

Let $Q=(y,z)k[y,z,\mathbf t]\cap R$. Then $QR[X]$ is a proper
t-ideal of $R[X]$ properly containing $P$. 

To verify this, let
$F$ be a finite subset of $QR[X]$. Let $t_n$ be an indeterminate
that does not occur in the polynomials in the set $F$. Then $t_n f\in
R[X]$ for all $f\in F$, so $t_n\in (R[X]:F)$. If $g\in (F)_v$, then
$gt_n\in R[X]$. Hence $\deg_{y,z}gt_n\ge 1$ and $g\in Q$. It
follows that $(F)_v\subseteq Q$, so $Q$ is a t-ideal.

\end{enumerate}

\end{proof}



\begin{ex} An example of a strong maximal divisorial ideal of an integrally 
closed domain $R$ that is not $t$-maximal.
\end{ex}

{\it Let $k$ be a field and let $Y,Z,\mathbf X=\left
\{X_n\,:\, n\ge1\right\},\mathbf T=\left\{T_n\,:\, n\ge1\right\}$ be
independent indeterminates over $k$.  Let $S$ be the
set of monomials $f$ in $k[Y,Z,\mathbf X,\mathbf T]$
satisfying the following two conditions:

\begin{enumerate}
\item[(a)]
If $Z$ occurs in $f$, then some $X_n$ occurs
in $f$.
\item[(b)]
For all $n$, if $T_n$ occurs in $f$, then either 
$Y$ or $X_i $ occurs in $f$ for some
$i\le n$.
\end{enumerate}

Clearly, $S$ is a semigroup containing $\mathbf X$ and $Y$.  Let $R=k[S]$ 
be the semigroup ring over $S$ and set
$$P=(\mathbf X)k[Y,Z,\mathbf X,\mathbf T]\cap R.$$ 
Then $R$ is integrally closed and $P$ is a strong maximal divisorial
ideal of $R$ that is not $t$-maximal.}

\begin{proof}
We will use repeatedly that $P$
is a monomial ideal of $R$.

\begin{enumerate}
\item
{\em $R$ is integrally closed.}

By \cite[Corollary 12.11 (2)]{gilmer}, it
is enough to show that the monoid $S$ is integrally closed.  If $f$ is
an element in the quotient group of $S$ such that $f^n\in S$ for some
$n\ge1$, then $f$ is a monomial.  Since $f^n$ satisfies conditions
(a)-(b), it is clear that $f$ also satisfies them, thus $f\in S$.
We conclude that $R$ is integrally closed.

\item
{\em 
$P=RZ^{-1}\cap R$. Hence $P$ is a divisorial 
ideal.}

 Clearly, any monomial in $ZP$ satisfies conditions
(a)-(b), hence $ZP\subseteq R$.  Thus $P\subseteq RZ^{-1}\cap R$.  

For the reverse inclusion, it is enough
to show that any monomial $f\in RZ^{-1}$ belongs to $P$. Since $Zf\in
R$, we see that $Zf$ satisfies conditions (a)-(b) and so does $f$,
thus $f\in R$.  Using again that $Zf\in R$ we see that some $X_n$
occurs in $f$, hence $f\in P$.

\item
{\em  $(R:P)=R[Z]$.}

Using conditions (a)-(b), we see that $R[Z]\subseteq
(R:P)$. 

For the reverse inclusion, let $u$ be a quotient of monomials in
$(R:P)$. Since $uX_1,uX_2\in R$, we see that $uX_1$ and $uX_2$ are
monomials, hence, by factoriality, $u$ also is a monomial. Let $u=Z^k
u_0$, where $k\ge0$, $u_0$ is a monomial and $Z$ does not occur in
$u_0$. Choose a positive integer $N$ such that $N>i$ for all $T_i$'s
occurring in $u$. Since $Z^ku_0X_N\in R$, we see that $u_0$ satisfies
condition (b); hence $u_0\in R$, so $u \in R[Z]$. 

\item
{\em  $P$ is a strong maximal divisorial ideal.}

We have $(R:P)=R[Z]\subseteq (P:P)$, thus
$(R:P)=(P:P)$, that is, $P$ is strong.

Assume that $P$ is not maximal divisorial, so there is a divisorial
ideal $Q$ properly containing $P$. Let $f\in Q\setminus P$. We may
assume that no $X_n$ occurs in $f$, thus $Z$ does 
not occur in $f$ either by condition (a) above.
Let $g\in (R:Q)\setminus R$, thus $g\in (R:P)=R[Z], g=\sum_{i=0}^n a_iZ^i$,
where $a_0,\dots, a_n\in R$. We may assume that $a_nZ^n\notin R$, thus 
$n\ge1$. We also may assume that
no $X_i$ occurs in $a_n$. Thus no $X_i$ occurs in $fa_n$,
which implies that $fa_nZ^n\notin R$. Since $R=k[S]$, we obtain
that $fg= fa_nZ^n+\dots\notin R$, a contradiction.

\item   {\em The ideal $M=(S)R$ is a maximal ideal of $R$ properly
	containing $P$ and is a t-ideal.}

	Clearly  $M$ is a maximal ideal containing $P$.
Since $Y\in M\setminus P$, we have $P\subsetneqq M$.

To show that $M$ is a $t$-ideal, let $F$ be a finite subset of $M$ and
let $N$ be a positive integer such that $N>i$ for each $T_i$ occurring
in some element of $F$. From conditions (a)-(b) it follows that
$M\subseteq (\mathbf X,Y)k[Y,Z,\mathbf X,\mathbf T]$. Hence
$T_NF\subseteq R$. Thus $(F)_v\subseteq (R:T_N)\cap R$. Since $T_N\notin
R$ and since $(R:T_N)\cap R$ is a monomial ideal, we obtain that
$(R:T_N)\cap R\subseteq (S)R=M$. It follows that $(F)_v \sub M$ and that
$M$ is a t-ideal. 

\end{enumerate}
\end{proof}




\end{document}